\documentclass[11pt]{article}

\usepackage{amssymb, enumerate}
\usepackage{amsfonts}
\usepackage{amsmath}
\usepackage{color}
\usepackage{hyperref}
\usepackage{epsfig}
\usepackage{amsthm}

\definecolor{green3}{rgb}{0,0.6,0}

\newtheorem{proposition}{Proposition}
\newtheorem{theorem}[proposition]{Theorem}
\newtheorem{definition}[proposition]{Definition}
\newtheorem{remark}[proposition]{Remark}

\newtheorem{notation}[proposition]{Notation}

\title{Other Examples of Principal Ideal Domains \\
that are not Euclidean Domains}
\author{Nicol\'as Allo-G\'omez$^{{*, \diamondsuit}}$}

\begin{document}

\maketitle

\noindent {\small ${*}$ Universidad de Buenos Aires. Facultad de Ciencias Exactas y Naturales. Departamento de Matem\'atica. Buenos Aires,  Argentina.}\\ 
{\small ${^\diamondsuit}$ Universidad Torcuato Di Tella. Departamento de Matemáticas y Estadística. Buenos Aires,  Argentina.}

\bigskip

\noindent E-mail: \texttt{nicolas.allo@utdt.edu}

\bigskip

\begin{abstract}
It is a well-known and easily established fact that every Euclidean domain is also a principal ideal domain. However, the converse statement is not true, and this is usually shown by exhibiting as a counterexample the ring of algebraic integers in a certain, very specific quadratic field, and the proof that this works is quite unnatural and technical. In this article, we will present a family of counterexamples constructed using real closed fields. 
\end{abstract}

\noindent \textbf{Keywords:} principal ideal domain, Euclidean domain, quotient of a polynomial ring, formally real field, real closed field.

\section{Introduction.}

Courses on modern algebra typically emphasize the importance of unique factorization domains (UFD for short) in solving equations. Two of the fundamental facts covered are that every principal ideal domain (PID) is a UFD, and that every Euclidean domain (that is, a domain with a division algorithm) is a PID. It is important to know, though, that the converse of these implications are not true. It is easy to produce examples of UFDs that are not PIDs. For instance, for any UFD $A$ that is not a field, the polynomial ring $A[X]$ is a UFD that is not a PID. However, it is more difficult to exhibit PIDs that are not Euclidean domains. There is a classical example that is most often mentioned: the ring of integers of the quadratic field $\mathbb{Q}(\sqrt{-19})$ (see, for example, \cite{Wilson, Campoli}). One should note that it is a very \emph{rigid} example in the sense that we cannot construct other examples using similar techniques, since it depends on certain arithmetical coincidences that occur in that ring. Another lesser-known example is the ring $\mathbb{R}[X,Y]/(X^2+Y^2+1)$. In \cite{samuel1}, P. Samuel proved that this ring is a UFD. More recently, in \cite{Stack}, there is a sketch of a proof that this ring is a PID but not a Euclidean domain, and a detailed proof of these facts can be found in \cite{pidnoted}, although it relies on certain non-elementary techniques and concepts. 

The aim of this paper is to provide a constructive and elementary proof of the essential result that is used to show that this example is indeed a PID and, by extension, to provide an immediate generalization that leads to a family of non-isomorphic examples of PIDs that are not Euclidean. With that goal, we first organize and present some necessary concepts and results, and then state and prove the main theorem of this paper, which is a description of the construction of \emph{other examples of PIDs that are not Euclidean domains}.

\section{Preliminaries.}

In this section, we present the concepts of \textit{formally real field} and \textit{real closed field}, which will play a fundamental role in our work. These serve as the cornerstone of modern semialgebraic geometry, which was formalized in the early 20th century by mathematicians E. Artin and O. Schreier in a long series of papers beginning with \cite{AS01, AS02}. We will also state a few of the classical results associated with these concepts, whose proofs, along with other complementary ones, can be found in \cite{BCR, BPR}.

Furthermore, we recall the notion of the resultant of two univariate polynomials, which will be an essential tool in the proof of our main theorem, and outline its most important properties.

\subsection{Real Fields.} In what follows we will consider fields endowed with order relations, we make the convention that these will always be total orders.

\begin{definition}
A field $F$ with an order $\leq$ is an ordered field if
\begin{itemize}
	\item $x \leq y  \Rightarrow x + z \leq y + z$,
	\item $x,y \geq 0 \Rightarrow x \cdot y \geq 0$
\end{itemize} 
for all $x,y,z \in F$. We say that a field $F$ can be ordered if there is an order $\leq$ in $F$ with respect to which $F$ is an ordered field.
\end{definition}

\begin{definition}
A formally real field $F$ is a field that satisfies one of the following equivalent conditions:
\begin{enumerate}
	\item $F$ can be ordered.
	\item The element $-1$ is not a sum of squares in $F$.
	\item If any sum of squares of elements of $F$ is zero, then each of those elements must be zero.
\end{enumerate} 
\end{definition}

Notice that, as a consequence of the third condition, $a^2 \neq -b^2$ whenever $a$ and $b$ are nonzero elements in a formally real field. We will use this below.

\medskip

Some examples of formally real fields are the field of rational numbers $\mathbb{Q}$, the real numbers~$\mathbb{R}$, and, more generally, any subfield $L$ of $\mathbb{R}$ such as the algebraic real numbers, and the hyperreal numbers (see \cite{Goldb}). Moreover, if $F$ is a formally real field, then $F(X)$, the field of rational functions over $F$, $F((X))$, the Laurent series over $F$, and $F\{\{X\}\}$, the Puiseux series over $F$, are all formally real fields.

\begin{definition}
A real closed field $F$ is a formally real field such that the extension $F(\sqrt{-1}) = F[X] / (X^2 + 1)$ is an algebraically closed field.
\end{definition}

Among the examples we listed above of formally real fields, the following are real closed fields:  $\mathbb{R}$, the algebraic real numbers, the hyperreal numbers and $\mathbb{R}\{\{X\}\}$. Besides that, one can show that if $F$ is a real closed field then the Puiseux series $F\{\{X\}\}$ over $F$ is also a real closed field, and, using the Artin-Schreier Theorem which can be found in \cite{Jac}, that every formally real field has a real extension that is real closed (see \cite[Theorem 2.2]{Lang}).

\begin{theorem}
The following conditions are equivalent for any field $F$:
\begin{enumerate}
	\item $F$ is a real closed field.
	\item $F$ is a formally real field that has no formally real proper algebraic extensions.
	\item $F$ can be ordered so that $\{a \in F : a \geq 0\}$ is the set of squares in $F$ and every polynomial $p \in F[X]$ of odd degree has at least one root in $F$. 
\qed
\end{enumerate}
\end{theorem}


\subsection{The resultant.} Let us turn now to the concept of the resultant of two univariate polynomials with coefficients in a UFD and some of its fundamental properties. The resultant dates back to the works of G.W. Leibniz, L. Euler, E. Bézout, and C.G.J. Jacobi, but its modern formulation is due to J.J. Sylvester in \cite{Sylv}. The resultant is an essential tool that arises naturally in many areas of mathematics and nowadays it also has significant algorithmic importance. It is a key tool of elimination theory, which was prominent until the mid-20th century and experienced a resurgence with the massive development of computational algebra in the late 1960s.

\begin{definition}
Let $A$ be a UFD and let $f$ and $g$ be nonconstant polynomials in $A[X]$ of degrees $n$ and $m$ respectively, so that $f = \displaystyle \sum_{i = 0}^n a_i X^i$, $g = \displaystyle \sum_{j = 0}^m b_j X^j$, with $a_n , b_m \neq 0$. The Sylvester matrix associated to $f$ and $g$ is the $(n+m) \times (n+m)$ matrix
{\small\[ 
\operatorname{Syl}(f,g) \ = 
\left[
\begin{array}{ccccccccccc}
a_n & a_{n-1} & \cdots & \cdots & \cdots & a_0 & 0 & \cdots & \cdots & \cdots & 0 \\
0 & a_n & a_{n-1} & \cdots & \cdots & \cdots & a_0& 0 & \cdots & \cdots & 0 \\
0 & 0 & a_n & a_{n-1} & \cdots & \cdots & \cdots & a_0 & 0 & \cdots & 0 \\
\vdots & \  & \ddots & \ddots & \ddots & \  & \ & \ & \ddots & \ddots & \vdots \\
\vdots & \  & \  & \ddots & \ddots & \ddots & \  & \ & \  & \ddots & 0 \\
0 & \cdots & \cdots & \cdots & 0 & a_n & a_{n-1} & \cdots & \cdots & \cdots & a_0 \\
b_m & b_{m-1} & \cdots & \cdots & \cdots & \cdots & b_0 & 0 & \cdots & \cdots & 0 \\
0 & b_m & b_{m-1} & \cdots & \cdots & \cdots & \cdots & b_0 & 0 & \cdots & 0 \\
\vdots & \ddots & \ddots & \ddots & \  & \  & \ & \ & \ddots & \ddots & \vdots \\
\vdots & \  & \ddots & \ddots & \ddots & \ & \ & \ & \  & \ddots & 0 \\
0 & \cdots & \cdots & 0 & b_m & b_{m-1} & \cdots & \cdots & \cdots & \cdots & b_0 
\end{array}
\right]
\]}
and the resultant of $f$ and $g$ is its determinant
$$\operatorname{Res}(f,g) = \det \Big(\operatorname{Syl}(f,g)\Big) \in A.$$
\end{definition}

We finish this section by presenting the two key properties of the resultant. Detailed proofs of these results, along with many others, can be found in \cite{VZGG, CLO}.

\begin{theorem}
Let $A$ be a UFD and let $f$ and $g$ be nonconstant polynomials in $A[X]$.
\begin{enumerate}
	\item $\operatorname{Res}(f,g) = 0$ if and only if there exists a nonconstant polynomial $h \in A[X]$ such that $h \mid f$ and $h \mid g$.
	\item There exist polynomials $s,t \in A[X]$ such that $\operatorname{Res}(f,g) = sf + tg$.
\qed
\end{enumerate}
\end{theorem}

\section{Principal Ideal Domains that are not Euclidean Domains.}
It is the aim of this article to prove the following main result:
\begin{theorem}\label{maintheo}
	Let $F$ be a real closed field. The ring $F[X,Y] / (X^2 + Y^2 + 1)$ is a PID that is not a Euclidean domain.
\end{theorem}

Let us fix the following notation for the rest of the paper.
\begin{notation}

Let $F$ be a real closed field. We set $A = F[X,Y] / (X^2 + Y^2 + 1)$ and we shall denote $\overline{f}$ as the class in $A$ of a polynomial $f \in F[X,Y]$. Besides, we will write $\operatorname{deg}_{X}(f)$ for the degree of $f$ as an element of $(F[Y])[X]$, and similarly for $\operatorname{deg}_{Y}(f)$.

For every commutative ring with unit $R$, we denote by $\mathcal{U}(R)$ its group of units.

\end{notation}

Note that $A$ is an integral domain since the polynomial $X^2 + Y^2 + 1$ is irreducible in $F[X,Y]$ and this ring is a UFD. Moreover, every $\overline{f} \in A$ can be written as $\overline{f} = \overline{q} \ \overline{X} + \overline{p}$ for a unique choice of $q, p \in F[Y]$, and since $A$ is a nontrivial $F$-algebra we may denote $\alpha$ the class in $A$ of an element $\alpha \in F$.

\medskip

We will now proceed to characterize the group of units of our ring $A$.

\begin{proposition}\label{unitsring}
	$\mathcal{U}(A) = F \setminus \{0\}$.
\end{proposition}

\begin{proof} Since every element $\alpha \in F \setminus \{0\}$ is clearly invertible, we need only focus on the other inclusion. Let $\overline{f}  =  \overline{q_1} \ \overline{X} + \overline{p_1}$ be an element in $\mathcal{U}(A)$ and let us prove that $q_1 = 0$ and $p_1 \in  F \setminus \{0\}$. Since $\overline{f}$ is invertible there is a $\overline{g} = \overline{q_2} \ \overline{X} + \overline{p_2} \in A$ such that $\overline{f} \ \overline{g} = 1$, which is equivalent to saying that the polynomial 
$$h = (q_1 q_2) \ X^2 + (q_1 p_2 + q_2 p_1) \ X + (p_1 p_2 - 1)$$
belongs to the ideal $(X^2 + Y^2 + 1)$ of $F[X,Y]$. Since $\operatorname{deg}_{X}(h) \leq 2$, we can write $h = r \ (X^2 + Y^2 + 1)$ for some $r \in F[Y]$, and thus we have the polynomial identities in $F[Y]$:
\begin{align}
q_1 q_2 & =  r, \\
q_1 p_2 + q_2 p_1& =  0, \\
p_1 p_2 & =  r \ (Y^2 + 1) + 1.
\end{align}
 
We will prove that $r$ must be the zero polynomial. Suppose that, contrary to our claim, $r \neq 0$. It then follows from identities (1) and (3) that $p_1, p_2, q_1, q_2 \neq 0$. Moreover, from identity (2) we have $q_1 p_2 = - q_2 p_1$.
Multiplying both sides of identity (3) by $q_1$ we obtain $q_1 p_1 p_2 = q_1 r \ (Y^2 + 1) + q_1$, and therefore $ - p_1^2 q_2 = q_1^2 q_2 \ (Y^2 + 1) + q_1$, which is equivalent to saying that 
$$q_2 \big(-p_1^2 - q_1^2 \ (Y^2 + 1)\big) = q_1.$$
Hence, $q_2$ divides $q_1$ in $F[Y]$.  Arguing in a similar way, by multiplying identity (3) by $q_2$ we can see that also $q_1$ divides $q_2$. But this means that there exists some $\lambda \in F \setminus \{0\}$ such that $q_2 = \lambda q_1$, and therefore $p_2 = - \lambda p_1$. From identity (3) we now get $- \lambda p_1^2 = \lambda q_1^2 \ (Y^2 + 1) + 1$, which is a contradiction: the leading coefficients of the polynomials in each side of this identity differ in a sign in the real closed field $F$. Therefore, $r = 0$ and the previous identities (1), (2) and (3) become
\begin{align*}
q_1 q_2 & =  0, \\
q_1 p_2 + q_2 p_1& =  0, \\
p_1 p_2 & =  1,
\end{align*}
which implies that $p_1, p_2 \in F \setminus \{0\}$ and $q_1 = q_2 = 0$, as desired.
\end{proof}

\medskip

In order to prove that the integral domain $A$ is a PID we are going to give a characterization of its nonzero prime ideals which, in fact, are maximal ideals. This characterization will also be useful later to conclude that this ring is not a Euclidean domain.

\begin{theorem}\label{primeideals}

A nonzero ideal $\overline{J}$ of $A$ is prime if and only if there exist $a,b \in F$ such that $\overline{J} = (\overline{X} + a\overline{Y} + b)$ \ or \ $\overline{J} = (\overline{Y} + a \overline{X} + b)$. Moreover, every nonzero prime ideal $\overline{J}$ of $A$ is maximal and $A / \overline{J} \simeq  F(\sqrt{-1})$.

\end{theorem}

\begin{proof} We will denote by $I$ the ideal of $F[X,Y]$ generated by $g = X^2 + Y^2 + 1$. 

Let us start by proving that every ideal of the form $\overline{J} = (\overline{X} + a\overline{Y} + b)$ with $a,b \in F$ is maximal and in particular prime. A similar reasoning applies to the other case.

Consider the ideal $J = (X + aY + b , X^2 + Y^2 + 1)$ of $F[X,Y]$. Clearly, $I \subseteq J$ and $\overline{J} = J/I$, so that we have the natural isomorphism 
$A/\overline{J} \ \simeq \ F[X,Y]/J$. Therefore, it suffices to show that $F[X,Y]/J$ is a field. 

Let $\overline{x}$ denote the image of $X$ in the quotient ring $F[X,Y]/J$, and similarly for $\overline{y}$. In this quotient we have $\overline{x} = - \left(a \overline{y} + b\right)$ and therefore 
\begin{equation*}
0 = (a^2+1) \overline{y}^2 + 2ab \overline{y} + (b^2 + 1).
\end{equation*}

Let now $p = (a^2+1) Y^2 + 2ab Y + (b^2 + 1) \in F[Y]$. Since $p \in J$, $p$ has degree $2$ and discriminant $\Delta(p) = -4 \left( a^2 + b^2 + 1 \right) < 0$ in $F$, we conclude that it has two different roots in $F(\sqrt{-1}) \setminus F$. 

Let $z \in F(\sqrt{-1}) \setminus F$ be one of these roots and consider the ring morphism \linebreak $\varphi : F[X,Y] \rightarrow F(\sqrt{-1})$ defined by $\varphi(f) = f\left(-az - b, z\right)$, which is surjective because it is $F-$linear and its image contains the element $z = \varphi(Y)$. We will prove that $\operatorname{Ker}(\varphi) = J$. Since $\varphi(X + aY + b) =  0$ and $\varphi(X^2 + Y^2 + 1) = p(z) = 0$, one inclusion is obvious. Let us deal with the other one. Let $f \in \operatorname{Ker}(\varphi)$ so that $f\left(-az - b, z\right) = 0$. We can write 
$$f = q \ \left(X + aY + b\right) + p s + r $$
for some $q \in F[X,Y], s, r \in F[Y]$, and either $r = 0$ or $\operatorname{deg}(r) < 2$. Since \linebreak $0 = \varphi(f) = r(z)$, it follows that $r = 0$ because $z \in F(\sqrt{-1}) \setminus F$ and $r$ has degree at most $1$. Therefore, $f = q \ \left(X + aY + b\right) + p s \in J$. 

Thus $J = \operatorname{Ker}(\varphi)$, as we said, and $\varphi$ induces an isomorphism $F[X,Y]/J  \simeq F(\sqrt{-1})$. This is a field, which is the desired conclusion.

\bigskip
We will next prove that every nonzero prime ideal of $A$ is as described in the theorem. Let $\overline{J}$ be such a nonzero prime ideal of $A$ and let us consider the prime ideal $J$ of $F[X,Y]$ such that $I \subsetneq J$ and $\overline{J} =  J/I$.
We will show that there is an element $q = aX + bY + c \in F[X,Y]$, with $(a,b) \neq (0,0)$, such that $q \in J$ and the assertion will follow from that. As before, let $\overline{x}$ and $\overline{y}$ denote the images of $X$ and $Y$, respectively, in the quotient ring $F[X,Y]/J$.

If either $X \in J$ or $Y \in J$, the result is clear. Let us suppose that instead $X, Y~\notin~J$, and thus that $\overline{x}, \overline{y} \neq 0$. Let $f \in J \setminus I$ and recall that $g = X^2 + Y^2 + 1$, so that $f\left(\overline{x},\overline{y}\right) = 0$ and $g\left(\overline{x},\overline{y}\right) = 0$ in $F[X,Y]/J$.

Since the quotient $F[X,Y]/J$ is an integral domain because $J$ is a prime ideal, we can consider its field of fractions $K =\operatorname{Frac}(F[X,Y]/J)$. Note that $F\subseteq F[X,Y]/J \subseteq K$, so we have a field extension $K / F$, and that the polynomial system
 \begin{equation} \tag{4} \label{systemres}
\left\{
\begin{aligned}
f(X,Y) & =  0, \\
g(X,Y) & =  0 \\
\end{aligned}
\right.
\end{equation}
has $\left(\overline{x}, \overline{y}\right)$ as a solution in $K^2$. We will prove that both $\overline{x}$ and $\overline{y}$ are algebraic over $F$.
Let us first consider the case in which $\operatorname{deg}_{X}(f) = 0$, so that $f \in F[Y]$ and $\operatorname{deg}_{Y} (f) > 0$. In this case it is clear that $\overline{y}$ is algebraic over $F$ because $f(\overline{y}) = 0$. Let us now take care of $\overline{x}$. Since $f \notin I = (g)$ and $g$ is irreducible in $F[X,Y]$, the polynomials $f$ and $g$ are coprime and therefore the resultant $R_Y = \operatorname{Res}(f,g,Y) \in F[X]$ of $f$ and $g$ viewed as polynomials in the variable $Y$ is nonzero. Recall that the resultant $R_Y$ satisfies B\'ezout's identity: we can write $R_Y = s f + t g$ for some polynomials $s,t \in F[X,Y]$. Since $\left(\overline{x}, \overline{y}\right)$ is a solution of the system (\ref{systemres}), we have that $R_Y (\overline{x}) = 0$, and thus that the element $\overline{x}$ is also algebraic over $F$. Similar arguments apply to the case in which $\operatorname{deg}_{Y}(f) = 0$, and the same conclusion can be drawn for the case where both $\deg_{X}(f), \deg_{Y}(f) > 0$ by considering the resultants $R_Y = \operatorname{Res}(f,g,Y) \in F[X]$ and $R_X = \operatorname{Res}(f,g,X) \in F[Y]$. 

Summarizing, we have proved that the elements $\overline{x}$ and $\overline{y}$ of $K$ are both algebraic over $F$, and thus that the field extension $F\big(\overline{x}, \overline{y} \big)/ F$ is algebraic. Since $F$ is a real closed field, the degree of that extension is at most $2$, and consequently the set $\big\{ \overline{x}, \overline{y}, 1 \big\}$ is linearly dependent over $F$. Therefore, there are $a,b,c \in F$, not all zero, such that $a \overline{x} + b \overline{y} + c = 0$, and this is equivalent to saying that $aX + bY + c$ \ belongs to $J$ . 
\end{proof}

\medskip

We can now proceed to show that our ring $A$ is in fact a PID. 

\begin{theorem}
The integral domain $A$ is a PID.
\end{theorem}

\begin{proof}\
Theorem \ref{primeideals} tells us that every prime ideal in $A$ is principal, and this implies that $A$ is a PID --- (see \cite{Kap, DuFoo}). 
\end{proof} 

From Theorem \ref{primeideals}, we also obtain the following characterization of the prime elements in $A$.

\begin{remark}\label{primelements}
Since $A$ is a PID in which all the nonzero prime ideals are of the form $(a\overline{X} + b\overline{Y} + c)$ with $(a,b) \neq (0,0)$, an element $\overline{p} \in A \setminus \{0\}$ is prime if and only if there exist $a,b,c \in F$, with $(a,b) \neq (0,0)$, such that $\overline{p} = a \overline{X} + b \overline{Y} + c$.
\end{remark}

In order to accomplish the aim of this paper, we are about to show that the PID under study is not a Euclidean domain. We will need first the following auxiliary result.

\begin{proposition}\label{unitsquot}
Let $R$ be a Euclidean domain that is not a field with a Euclidean function $\varphi: R \setminus \{0\} \longrightarrow \mathbb{Z}_{\geq 0}$. There exists a prime element $p \in R$ such that $\pi \big(\mathcal{U}(R)\big) = \mathcal{U}\big(R/(p)\big)$, where $\pi : R \longrightarrow R/ (p)$ is the natural projection.
\end{proposition}

\begin{proof}

Even though this result is well established in the literature (see \cite{samuel2}), we include a proof here for the sake of completeness.

Since $R$ is not a field, we can choose a nonzero element $p \in R$ that is not a unit and such that $\varphi(p)$ is minimal. We will show that $p$ is necessarily prime and that it has the desired property.


Let $a,b \in R$ such that $p$ divides $ab$, but $p$ does not divide $a$. Then there exist $q, r~\in~R$, $r \neq 0$, such that $a = qp + r$, with $\varphi(r) < \varphi(p)$. By the minimality of $\varphi(p)$, it follows that $r$ is a unit in $R$. Since $ab = qpb + rb$ and $p$ divides $ab$, we conclude that $p$ divides $rb$. As $r$ is a unit, this implies that $p$ divides $b$.

Now consider the natural projection $\pi : R \longrightarrow R/ (p)$. We will show that $\pi \big(\mathcal{U}(R)\big)$ $= \mathcal{U}\big(R/(p)\big)$.
Let $\overline{u}$ be an element in $\mathcal{U}\big(R/(p)\big)$ and let us consider $u \in R$ such that $\pi(u) = \overline{u}$. Since $\overline{u} \neq 0$ in $R/(p)$, it follows that $p$ does not divide $u$ in $R$. Hence, there exist $q, r \in R$, with $r \neq 0$, such that $u = qp + r$ with $\varphi(r) < \varphi(p)$. Therefore, $\overline{u} = \pi(r)$, and by the minimality of $\varphi(p)$, we conclude that $r$ is a unit in $R$. The other inclusion is straightforward.

\end{proof}

We can finally prove that $A$ is not a Euclidean domain, and with this reach the goal of this paper. 

\begin{theorem}

The integral domain $A$ is not a Euclidean domain.

\end{theorem}

\begin{proof} 
Assume $A$ is a Euclidean domain. Since $A$ is not a field, it follows from Proposition~\ref{unitsquot} and the characterization given in Remark~\ref{primelements} that there exists a prime element $\overline{p}$ of the form $a \overline{X} + b \overline{Y} + c \in A$, with $(a,b) \neq (0,0)$, such that $\pi \big(\mathcal{U}(A)\big) = \mathcal{U}\big(A/(\overline{p})\big)$, where $\pi : A \longrightarrow A/(\overline{p})$ is the natural projection. Restricting $\pi$ yields a surjective group homomorphism $\psi : \mathcal{U}(A) \longrightarrow \mathcal{U}\big(A/(\overline{p})\big)$. We will prove that $\psi$ is also injective. 

On account of Proposition \ref{unitsring} we have that $\mathcal{U}(A) = F \setminus \{ 0 \}$, and this implies that $A$ has the unusual property that whenever $\alpha$ and $\beta$ are different elements of $\mathcal{U}(A)$ then $\alpha - \beta \in \mathcal{U}(A)$.
Let $\alpha, \beta \in \mathcal{U}(A)$ be such that $\alpha \neq \beta$ and suppose $\psi(\alpha) = \psi(\beta)$, so $\overline{p} \mid \alpha - \beta$ in $A$. The element $\overline{p}$ is prime in $A$ and $\alpha - \beta \in \mathcal{U}(A)$, so this is impossible and therefore $\psi$ is injective. 

As a consequence of all this, we have a group isomorphism $F \setminus \{ 0 \} \simeq \mathcal{U}\big(A/(\overline{p})\big)$. Since the rings $A/ (\overline{p})$ and $F(\sqrt{-1})$ are isomorphic, it follows that $\mathcal{U}\big(A/(\overline{p})\big)$ and $F(\sqrt{-1}) \setminus \{ 0 \}$ are isomorphic as groups, so $F \setminus \{ 0 \}$ and $F(\sqrt{-1}) \setminus \{ 0 \}$ are also isomorphic. This implies that there is an element $r \in F$ such that $r^2 = -1$, and this contradicts the fact that $F$ is a real closed field. 
This contradiction arose from our assumption that $A$ is a Euclidean domain, so it must not be one.
\end{proof}

As we noted at the beginning of this article, we can use Theorem \ref{maintheo} to produce examples of non-isomorphic PIDs that are not Euclidean. If we start with a real closed field $F$ then the subring of $A$ generated by the units of $A$ is $F$ itself. Therefore, if two real closed fields $F$ and $F'$ are not isomorphic then the corresponding rings $A$ and $A'$ are themselves not isomorphic.

\medskip
\bigskip
\bigskip

\noindent \textbf{Acknowledgments.}
The author is deeply indebted to Teresa Krick, Mariano Su\'arez-\'Alvarez, and Juan Sabia for their many helpful suggestions and for their active interest in the publication of this paper. The author also wishes to thank the anonymous reviewers for their careful reading of the manuscript and their valuable comments and suggestions.

\medskip

\bibliographystyle{vancouver}
\bibliography{PIDnotEDarxi}

\begin{thebibliography}{10}

\bibitem{Wilson}
Wilson JC.
\newblock A principal ideal ring that is not a {E}uclidean ring.
\newblock Math Mag. 1973;46:34-8.
\newblock Available from: \url{https://doi.org/10.2307/2688577}.

\bibitem{Campoli}
C\'ampoli OA.
\newblock A principal ideal domain that is not a {E}uclidean domain.
\newblock Amer Math Monthly. 1988;95(9):868-71.
\newblock Available from: \url{https://doi.org/10.2307/2322908}.

\bibitem{samuel1}
Samuel P.
\newblock Anneaux Factoriels.
\newblock Sociedade de Matem\'atica de S\~ao Paulo; 1963.

\bibitem{Stack}
user26857.
\newblock Quotient of polynomials, PID but not Euclidean domain? [webpage];
  2017.
\newblock Available from: \url{https://math.stackexchange.com/q/864627}.
\newblock URL:https://math.stackexchange.com/q/864627 (version: 2017-04-13).

\bibitem{pidnoted}
Class-notes.
\newblock Non-Euclidean PID [webpage]; 2016.
\newblock Available from: \url{https://www.mat.uniroma2.it/~
  eal/NoneuclPID.pdf}.
\newblock URL: https://www.mat.uniroma2.it/~ eal/NoneuclPID.pdf.

\bibitem{AS01}
Artin E, Schreier O.
\newblock Algebraische {K}onstruktion reeller {K}\"orper.
\newblock Abh Math Sem Univ Hamburg. 1927;5(1):85-99.
\newblock Available from: \url{https://doi.org/10.1007/BF02952512}.

\bibitem{AS02}
Artin E, Schreier O.
\newblock Eine {K}ennzeichnung der reell abgeschlossenen {K}\"orper.
\newblock Abh Math Sem Univ Hamburg. 1927;5(1):225-31.
\newblock Available from: \url{https://doi.org/10.1007/BF02952522}.

\bibitem{BCR}
Bochnak J, Coste M, Roy MFc.
\newblock Real algebraic geometry. vol.~36 of Ergebnisse der Mathematik und
  ihrer Grenzgebiete (3).
\newblock Springer-Verlag, Berlin; 1998.
\newblock Translated from the 1987 French original, Revised by the authors.
\newblock Available from: \url{https://doi.org/10.1007/978-3-662-03718-8}.

\bibitem{BPR}
Basu S, Pollack R, Roy MFc.
\newblock Algorithms in real algebraic geometry. vol.~10 of Algorithms and
  Computation in Mathematics.
\newblock 2nd ed. Springer-Verlag, Berlin; 2006.

\bibitem{Goldb}
Goldblatt R.
\newblock Lectures on the hyperreals: an introduction to nonstandard analysis.
\newblock Springer; 1998.

\bibitem{Lang}
Lang S.
\newblock XI.
\newblock In: Algebra. third (revised) ed. Springer-Verlag; 2002. p. 451-3.

\bibitem{Jac}
Jacobson N.
\newblock 11.
\newblock In: Basic Algebra II. W. H. Freeman and Co.; 1989. p. 650-7.

\bibitem{Sylv}
Sylvester JJ.
\newblock On a theory of the syzygetic relations of two rational integral
  functions, comprising an application to the theory of Sturm’s functions,
  and that of the greatest algebraical common measure, Philosophical
  Transactions of the Royal Society of London, Part III (1853), pp. 407– 548.
\newblock In: Collected Mathematical Papers of James Joseph Sylvester. vol.~1.
  Chelsea Publishing Co.; 1973. p. 429-586.

\bibitem{VZGG}
von~zur Gathen J, Gerhard J.
\newblock Modern computer algebra.
\newblock 3rd ed. Cambridge University Press, Cambridge; 2013.
\newblock Available from: \url{https://doi.org/10.1017/CBO9781139856065}.

\bibitem{CLO}
Cox DA, Little J, O'Shea D.
\newblock Ideals, varieties, and algorithms.
\newblock 4th ed. Undergraduate Texts in Mathematics. Springer, Cham; 2015.
\newblock An introduction to computational algebraic geometry and commutative
  algebra.
\newblock Available from: \url{https://doi.org/10.1007/978-3-319-16721-3}.

\bibitem{Kap}
Kaplansky I.
\newblock Elementary divisors and modules.
\newblock Trans Amer Math Soc. 1949;66:464-91.
\newblock Available from:
  \url{https://doi.org/10.1090/S0002-9947-1949-0031470-3}.

\bibitem{DuFoo}
Dummit DS, Foote RM.
\newblock Abstract Algebra.
\newblock 3rd ed. John Wiley $\&$ Sons; 2004.
\newblock Section 8.2, Exercise 10.

\bibitem{samuel2}
Samuel P.
\newblock About Euclidean rings.
\newblock Journal of Algebra. 1971;19:282-301.

\end{thebibliography}

\end{document}